\documentclass{article}
\usepackage[utf8]{inputenc}
\usepackage[margin=1in]{geometry}
\usepackage{amsmath} 
\usepackage{amsthm}
\usepackage{todonotes}
\usepackage[nolist]{acronym}
\usepackage{xcolor}
\usepackage{amsmath}
\usepackage{amssymb}
\usepackage{mathtools}
\usepackage{subcaption}
\usepackage{hyperref}
\hypersetup{colorlinks=true, linkcolor=black}
\usepackage[ruled]{algorithm2e}
\usepackage{mathrsfs}

\newcommand{\R}{\mathbb{R}}

\newcommand{\N}{\mathbb{N}}

\newcommand{\0}{\mathbf{0}}
\newcommand{\1}{\mathbf{1}}

\newcommand{\dc}{\mathcal{D}}

\DeclarePairedDelimiter{\abs}{\lvert}{\rvert}
\DeclarePairedDelimiter{\norm}{\lVert}{\rVert}

\DeclarePairedDelimiter{\diag}{\textrm{diag}(}{)}
\DeclarePairedDelimiter{\vvec}{\textrm{vec}(}{)}

\DeclarePairedDelimiterX{\inp}[2]{\langle}{\rangle}{#1, #2}

\DeclareMathOperator*{\find}{find}

\newtheorem{thm}{Theorem}[section]

\newtheorem{lem}[thm]{Lemma}

\newtheorem{prob}[thm]{Problem}

\newtheorem{defn}{Definition}[section]

\newtheorem{rmk}{Remark} 
\newcommand{\bx}{\mathbf{X}}
\newcommand{\bu}{\mathbf{U}}
\newcommand{\bw}{\mathbf{W}}
\newcommand{\bs}{\mathbf{S}}

\newcommand{\xn}{\mathbf{X}_-}
\newcommand{\xd}{\mathbf{X}_\delta}
\newcommand{\bth}{\mathbf{\Theta}}
\newcommand{\dx}{\delta x}

\newcommand{\scr}{\mathcal{S}}

\newcommand{\kron}{\otimes}
\newcommand{\kcol}{\otimes_{\textrm{col}}}

\title{\LARGE \bf Data-Driven Control of Positive Linear Systems \\ using Linear Programming
}

\renewcommand\footnotemark{}

\author{Jared Miller $^1$, Tianyu Dai $^2$,  Mario Sznaier $^1$, Bahram Shafai $^1$
\thanks{$^1$J. Miller and M. Sznaier are with the Robust Systems Lab,  ECE Department, Northeastern University, Boston, MA 02115. (e-mails: miller.jare@northeastern.edu, msznaier@coe.neu.edu). B. Shafai is with the ECE Department, Northeastern University (e-mail: shafai@coe.neu.edu).}
\thanks{$^2$ T. Dai is with The MathWorks, Inc., 1 Apple Hill Drive,
		Natick, MA 01760 USA (e-mail: tdai@mathworks.com)}
\thanks{J. Miller and M. Sznaier were partially supported by NSF grants  CNS--1646121, ECCS--1808381 and CNS--2038493, AFOSR grant FA9550-19-1-0005, and ONR grant N00014-21-1-2431.  
}
}

\begin{document}

\maketitle


\begin{abstract}
\label{sec:abstract}
This paper presents a linear-programming based algorithm to perform data-driven stabilizing control of linear positive systems. A set of state-input-transition observations is collected up to magnitude-bounded noise. A state feedback controller and dual linear copositive Lyapunov function are created such that the set of all data-consistent plants is contained within the set of all stabilized systems. This containment is certified through the use of the Extended Farkas Lemma and solved via Linear Programming. Sign patterns and sparsity structure for the controller may be imposed using linear constraints. The complexity of this algorithm scales in a polynomial manner with the number of states and inputs. The Linear-Programming based algorithm is extended to positive-stabilization of switched linear  systems and Linear Parameter-Varying systems. Effectiveness is demonstrated on example systems.

\end{abstract}
\section{Introduction}
\label{sec:introduction}

This paper performs \ac{DDC} of Positive \ac{LTI} \acp{CTS} and \acp{DTS} by finding full-state-feedback stabilizing controllers. These controllers, which stabilize all possible plants that are consistent with observed data, are formulated as the solution to \iac{LP}. 

Positive systems are a class of dynamical systems whose state and output responses to positive (nonnegative) initial conditions and inputs remain positive (nonnegative) 
 for all time \cite{luenberger1979dynamic, berman1989nonnegative, farina2000positive, kaczorek2012positive}. Instances of positive systems include population models \cite{hirsch1988systems}, chemical networks \cite{blanchini2014piecewise}, radio communications \cite{ZAPPAVIGNA2012219}, queuing \cite{shorten2006positive}, 
and Markov chains \cite{seneta2006non}. 
Full-state-feedback stabilization of known \ac{LTI} positive systems can be accomplished by solving \iac{LP} to find control  (dual) linear copositive Lyapunov functions \cite{rami2007controller}.
Alternatively, one can perform stabilization by formulating \iac{SDP} to find a quadratic Lyapunov function \cite{gao2005control, shafai2013stabilityradius}.

The peak-to-peak ($L_\infty \rightarrow L_\infty$ for \iac{CTS} or $\ell_\infty \rightarrow \ell_\infty$ for \iac{DTS}) gain of an extended positive plant can be calculated and regulated using \iac{LP} \cite{ebihara2011l1, briat2013positive, ebihara2013stability, naghnaeian2014linfinity}, which has also been derived using stability radius formulas \cite{shafai2019stability}.
Analysis and stabilization results can be extended to uncertain and switched positive systems \cite{blanchini2015switched}, as well as time-delay positive systems \cite{shafai2014positive}. The tutorial in \cite{rantzer2018tutorial} is a  survey of topics about stabilization and performance regulation for positive linear systems.

\ac{DDC} is a method that synthesizes controllers for a class of data-consistent plants without first performing a possibly expensive and inaccurate system identification step \cite{HOU20133, hou2017datasurvey}. 
Methods that require a reference signal include iterative feedback tuning \cite{hjalmarsson1998iterative}, virtual reference feedback tuning \cite{campi2002virtual}, \cite{bazanella2011data}, and correlation-based tuning \cite{karimi2004iterative}, but these algorithms lack stability guarantees for all consistent systems.
Data-driven predictive control through input-output data can be accomplished through Willem's Fundamental Lemma, assuming that a rank condition of the Hankel matrices is satisfied (persistency of excitation) \cite{willems2005note}. Stabilization,  worst-case-optimal control, and Model Predictive Control problems can be solved through the use of this Lemma \cite{waarde2020informativity, depersis2020formulas, coulson2019data, berberich2020mpc}, but the Lemma is vulnerable to noise sensitivity (even with regularization).

Prior knowledge of noise characteristics can be employed to synthesize controllers that will stabilize all plants that are consistent with data. $L_\infty$-bounded noise arises from  bounds on the time-derivative of the state (\ac{CTS}) or discretization of continuous-time finite-difference approximations (\ac{DTS}). Work addressing \ac{DDC} of $L_\infty$-bounded noise by solving \acp{LP} includes \cite{cheng2015} using an Extended Farkas Lemma \cite{hennet1989farkas}. Tools from polynomial optimization may be applied to the $L_\infty$ setting, such as for quadratic stabilization \cite{dai2020data}, switched systems \cite{dai2018moments, dai2018data}, and error-in-variables control \cite{miller2022eiv_short, miller2022eivarx}.
Quadratic Matrix Inequalities may be used to represent consistency sets (including energy-based or $L_2$-bounded noise) \cite{waarde2020noisy, berberich2020combining}, and stabilizing controllers may be synthesized by solving \acp{SDP} using a Matrix S-Lemma \cite{yakubovich1997s}. The work in \cite{martin2021data} employs polynomial optimization for \ac{DDC} under the assumption that magnitude bounds on Taylor polynomial coefficients and residual terms are known.

The work in \cite{shafai2022data} utilizes the Fundamental Lemma \cite{willems2005note} to perform \ac{DDC} of positive systems by solving \iac{SDP}. System identification of positive systems is performed in \cite{grussler2017identification}. The method in \cite{bianchi2022data} uses data-driven Lyapunov-Metzler inequalities to perform switched positive-systems control at the expense of solving Bilinear Matrix Inequalities.






The contributions of this work are:
\begin{itemize}
    \item \Iac{LP} that performs data-driven positive-stabilizing control for all systems consistent with observed data.
    \item A tabulation of computational complexity required to solve this \ac{LP}.
    \item An extension of this \ac{LP} towards worst-case peak-to-peak gain minimizing control.
    \item An application of the \ac{LP} for positive switched and positive \ac{LPV} systems.
\end{itemize}

This paper has the following structure: 
Section \ref{sec:preliminaries} reviews the preliminaries of notation, positive systems, copositive Lyapunov functions, and the Extended Farkas Lemma.
Section \ref{sec:stab} presents \iac{LP} to perform data-driven stabilizing control of positive systems. Section \ref{sec:peaktopeak} extends this \ac{LP} framework to yield controllers that minimize the worst-case peak-to-peak gain between an external input and a controlled output. Section \ref{sec:switched} extends the \ac{LP} stabilization methods towards control of switched positive linear systems. Section \ref{sec:lpv} details how these approaches can be used for positive \ac{LPV} systems.
Section \ref{sec:examples} demonstrates effectiveness of these methods on stabilizing and worst-case-optimal control of example systems.
Section \ref{sec:conclusion} concludes the paper.
\section{Preliminaries}
\label{sec:preliminaries}

\begin{acronym}[DLCLF]

\acro{CTS}{Continuous-Time System}

\acro{DDC}{Data-Driven Control}

\acro{DTS}{Discrete-Time System}





\acro{DLCLF}{Dual Linear Copositive Lyapunov Function}
\acro{LCLF}{Linear Copositive Lyapunov Function}
\acroindefinite{LCLF}{an}{a}

\acro{LP}{Linear Program}
\acroindefinite{LP}{an}{a}

\acro{LPV}{Linear Parameter-Varying}
\acroindefinite{LPV}{an}{a}

\acro{LPVA}{Linear Parameter-Varying A-Affine}
\acroindefinite{LPVA}{an}{a}

\acro{LTI}{Linear Time Invariant}
\acroindefinite{LTI}{an}{a}







\acro{SDP}{Semidefinite Program}
\acroindefinite{SDP}{an}{a}



\end{acronym}

\subsection{Notation}

The $n$-dimensional Euclidean vector space is $\R^n$. Its nonnegative orthant will be written as $\R^n_{\geq 0}$ and its positive orthant will be denoted as $\R^n_{> 0}$. The set of $n \times m$ matrices will be $\R^{n \times m}$. The transpose of a matrix $M \in \R^{n \times m}$ is $M^T \in \R^{m \times n}$.

The $n$-dimensional identity matrix is $I_n$. The vector of all ones is $\1_n \in \R^n$. The $n \times m$ matrix of all zeros is $\0_{n \times m} \in \R^{n \times m}$.
The matrix with $v \in \R^n$ appearing on its main diagonal and zeros elsewhere is $\diag{v} \in \R^{n \times n}$.
The Kronecker product of matrices $A$ and $B$ is $A \otimes B$. The column-wise vectorization of a matrix $M$ is $\vvec{M}$.  The elementwise division between $a, b \in \R^n$ is $a./b$.

The symbol $\delta x$ will refer to $x^+$ (next state) in discrete-time or $\dot{x}$ in continuous-time.
The symbols $(\circledast, \oplus, \ominus, \odot)$ correspond to an unrestricted (real-valued), a nonnegative, a nonpositive, and a zero-valued element respectively.

\subsection{Positive Systems}

A controlled \ac{LTI} system with states $x \in \R^n$, inputs $u \in \R^m$, and outputs $y\in\R^p$ has the form
\begin{align}
\label{eq:sys}
    \delta x &= A x + B u & y&=C x + D u.
\end{align}

\subsubsection{Positive System Descriptors}

\begin{defn}
The system \eqref{eq:sys} is \textit{internally positive} iff for any initial condition $x(0) \in \R^n_{\geq 0}$ and input $u(t) \in \R^m_{\geq 0}$, the state and output responses remain in the positive orthant ($x(t) \in \R^n_{\geq 0}, \ y(t) \in \R^p_{\geq 0} \ \forall t \geq 0$) \cite{farina2000positive}.
\end{defn}

Internal positivity requires that $(B, C, D)$ are all nonnegative, along with the property that $A$ is Metzler (off-diagonals are nonnegative) for \iac{CTS} or that $A$ is  nonnegative for \iac{DTS}. The system is positive-stable if $A$ is Hurwitz and Metzler (\ac{CTS}), or Schur and Nonnegative (\ac{DTS}).
For the remainder of this paper, we will assume that $C = I_n$ and $D = \0_{n \times m}$.




The state-feedback control $u= K x$ with $K \in \R^{m \times n}$ positively-stabilizes \eqref{eq:sys} if the closed-loop matrix $A+BK$ is Metzler-Hurwitz or Nonnegative-Schur (as appropriate).



\subsubsection{Copositive Functions}

\begin{defn}
    A function $f: \R^n \rightarrow \R$ is \textit{copositive} (with respect to the positive orthant) if $\forall x \in \R^n_{> 0}: f(x) > 0$. 
    \end{defn}
    
    Copositivity of the linear function $V(x)=v^T x$ and the dual-linear function  $V(x) = \max(x./v)$ may be checked by verifying that $v > 0$, but copositivity of a matrix function such as $x^T M x$ for some $M \in \R^{n \times n}$ is generically NP-hard \cite{murty1987some}.

\subsubsection{Stability of Positive Systems}
\begin{thm}
\label{thm:stability_pos}
    Let the system \eqref{eq:sys} be internally positive. Then it is asymptotically stable (positive-stable) iff one of the following equivalent conditions is satisfied \cite{rantzer2018tutorial}:
    \begin{enumerate}
        \item[C1)] The matrices $-A$ (\ac{CTS}) or $I_n - A$ (\ac{DTS}) have positive principal minors.
        \item[C2)] There exists a $p \in \R^n_{>0}$ with $P=\diag{p}$ such that $A^T P + P A \prec 0$ (\ac{CTS}) or $A^T P A - P \prec 0$  (\ac{DTS}).
        \item[C3)] There exists a positive vector $v \in \R^n_{>0}$ with \iac{LCLF} $v^T x$ such that $A^T v  < 0$ (\ac{CTS}) or $A^T v  < v_1$ (\ac{DTS}).        
        \item[C4)] There exists a positive vector $v_\infty \in \R^n_{>0}$ with \iac{DLCLF} $\max(x./v_\infty)$ such that $A v_\infty  < 0$ (\ac{CTS}) or $A v_\infty  < v_\infty$ (\ac{DTS}).
    \end{enumerate}
\end{thm}

In this paper we will exclusively use Condition C4 of Theorem \ref{thm:stability_pos} with \iac{DLCLF} $\max(x./v_\infty)$. We note that the conditions in Theorem \ref{thm:stability_pos} strictly treat the case of (dual) \acp{LCLF}. Proposition 3.3 of \cite{blanchini2015switched} states that every uniformly exponentially stable positive linear system admits a polyhedral Lyapunov function with an undecidable number of facets.


\subsubsection{Positive System Stabilization}

\acp{DLCLF} may be employed to find positive-stabilizing controllers $K \in \R^{m \times n}$. 
\begin{thm}[\cite{rami2007controller}]
\label{thm:stab_control}
The closed-loop system $\delta x = (A + BK) x$ from \eqref{eq:sys}, given a control $u=Kx$, is positive and asymptotically stable if there exists a vector $v \in \R^n_{>0}$ with a diagonal matrix $X = \diag{v}$, and a matrix $Y \in \R^{m \times n}$ such that the gain $K$ satisfies $K X = Y $ and
\begin{subequations}
\label{eq:stab_clean}
\begin{align}
    -(A X + B Y) \1_n & \in \R^n_{>0}  & &AX+BY  \text{ is Metzler (\ac{CTS})} \label{eq:stab_clean_c}\\
    v-(A X + B Y) \1_n & \in \R^n_{>0}  & & AX+BY \in \R^{n \times n}_{\geq 0} \ \text{(\ac{DTS})}. \label{eq:stab_clean_d}
\end{align}
\end{subequations}
\end{thm}



Finding a controller through \eqref{eq:stab_clean} requires solving \iac{LP} with both strict and nonstrict inequality constraints.

\subsubsection{Structured Control}

The stabilization task in \eqref{eq:stab_clean} may be restricted to a set of controllers that obey sign patterns and sparsity structures. Such sparsity might arise from network information constraints.

Let $\scr$ be an $m \times n$ matrix filled with the symbols $(\circledast, \oplus, \ominus, \odot)$. A controller with the structure $K \in \scr$ may be constructed by solving \eqref{eq:stab_clean} under the constraint that $Y \in \scr$, given that multiplication by the matrix $X$ with $v \in \R_{>0}^n$ does not change the sign pattern. An unstable internally positive system cannot be positive-stabilized by a nonnegative state feedback controller $K \in \R^{m \times n}_{\geq 0}$.

\subsection{Extended Farkas Lemma}

This work will find a state-feedback controller $u=Kx$ such that the set of all $K$-stabilized systems contains the set of systems consistent with observed data. The method used to enforce this containment is the Extended Farkas Lemma: 

\begin{lem}[Extended Farkas Lemma \cite{hennet1989farkas, henrion1999control}]
\label{lem:ext_farkas}
Let $P_1 = \{x \mid G_1 x \leq h_1\}$ and $P_2 = \{x \mid G_2 x \leq h_2\}$ be a pair of polytopes. Then $P_1 \subseteq P_2$ if and only if there exists a nonnegative matrix $Z$ of compatible dimensions such that,
\begin{align}
    Z G_1 &= G_2, & Z h_1 \leq h_2.
\end{align}
\end{lem}




\section{Data-Driven Stabilization}
\label{sec:stab}

This section will detail the data-driven positive-stabilization problem and its solution using robust linear programming.

\subsection{Problem Setting}

A set of $T$ observations are recorded of system \eqref{eq:sys} as corrupted by a noise process $w \in \R^n$,
\begin{equation}
    \delta x(t) = A x(t) + B u(t) + w(t).\label{eq:sys_corrupt}
\end{equation}

These observations are collected into the data $\dc = (\bx, \bu, \xd)$ with the expressions,
\begin{align}
\label{eq:data}
    \begin{array}{ccllll}
        \bx & := & [x(0) & x(1) &  \ldots & x(T-1) ]   \\
        \bu & := & [u(0) & u(1) & \ldots & u(T-1)] \\ 
        \xd & := & [\delta x(0) & \delta x(1) &  \ldots & \delta x(T-1) ].
    \end{array}
\end{align}

The discrepancy matrix $\bw$ satisfies the relation,
\begin{equation}
\label{eq:discrepancy}
    \bw = \xd - (A \bx + B \bu).
\end{equation}

The noise model that we will use is that each $w(t)$ (column of $\bw$) is $L_\infty$-norm-bounded by some given $\epsilon \geq 0$ ($\norm{w(t)}_\infty \leq \epsilon)$.

The set of all system matrices $(A, B)$ that are compatible with the $L_\infty$-corrupted data in $\dc$ forms a polytopic consistency set $\Sigma_{\dc}$. 
If it is known \textit{a priori} that $A$ is Metzler/Nonnegative and/or $B$ is nonnegative, then these constraints in $(A, B)$ may be adjoined to $\Sigma_\dc$.

The data-driven positive-stabilization problem is:
\begin{prob}
\label{prob:stab}
Find a vector $v \in \R^n_{>0}$ and a controller $K \in \scr$ such that $\max(x./v)$ is a common \ac{DLCLF} ensuring positive-stability of 
$A+BK$ for all $(A, B) \in \Sigma_\dc$.
\end{prob}


\subsection{Polytope Description}
We will describe $K$-stabilized and $\dc$-consistent polytopes that will be used in solving Problem \ref{prob:stab}
Throughout this section, the column-vectorization of the plant matrices will be defined as $a = \vvec{A}, b = \vvec{B}$. The identity $\vvec{U V W} = (W^T \otimes U) \vvec{V}$ for matrices $(U, V, W)$ of compatible dimensions will be judiciously used in derivations.

\subsubsection{Data-Consistent Polytopes}

The polytopic set $\Sigma_\dc^{\text{data}}$ of plants consistent with the data in $\dc$ may be represented as
\begin{subequations}
\label{eq:data_single}
\begin{align}
    G^{\textrm{data}}_1 &= \begin{bmatrix} \bx^T \otimes I_n & \bu^T \otimes I_n \end{bmatrix}\\
    \Sigma_\dc^{\text{data}} &= \left\{ (A, B) \mid G^{\textrm{data}}_1 \begin{bmatrix} a \\ b\end{bmatrix} \leq   \begin{bmatrix}\epsilon\1_{nT} +\vvec{\xd} \\ \epsilon\1_{nT} -\vvec{\xd}\end{bmatrix} \right\}.
\end{align}
\end{subequations}

The consistency set of plants $\Sigma_\dc$ is the intersection of $\Sigma_\dc^{\textrm{data}}$ and the prior knowledge on system-positivity of $(A, B)$ (linear constraints) described in $\Sigma^{\text{prior}}$. As an example, where $A$ is a positive system in discrete-time, then $\Sigma^{\text{prior}} = \{A \mid A \in \R^{n \times n}_+\}$, $G_1^{\textrm{prior}} = -I_{n^2}$, and $h_1^{\textrm{prior}} = \0_{n^2}$.
Let $(G_1, h_1)$ be matrices such that the polytopic data-consistency set $\Sigma_\dc = \Sigma^{\textrm{data}}_\dc \cap \Sigma^{\textrm{prior}}$ can be expressed as
\begin{subequations}
\begin{align}
    \Sigma_\dc = P_1 =  \left\{ (A, B) \mid G_1 \begin{bmatrix} a \\ b\end{bmatrix} \leq h_1 \right\}. \label{eq:sigma_dc}
\end{align}
\end{subequations}

\subsubsection{Controller-Stabilizing Polytopes}
\label{sec:stab_poly}
In order to apply the Extended Farkas Lemma \ref{lem:ext_farkas}, we will convert the strict inequalities in \eqref{eq:stab_clean} and in $v \in \R^n_{>0}$ to non-strict inequalities by utilizing a sufficiently small $\eta > 0$.
\begin{subequations}
\label{eq:stab_delta}
\begin{align}
    -(A X + B Y) \1_n -\eta \1_n& \in \R_{\geq 0}  & & \text{(\ac{CTS})} \label{eq:stab_delta_c} \\
    v-(A X + B Y) \1_n -\eta \1_n& \in \R_{\geq 0} & & \text{(\ac{DTS})} \label{eq:stab_delta_d}.
\end{align}
\end{subequations}

Define the canonical Metzler-indexing matrix $M_n \in \R^{n(n-1) \times n^2}$ as a 0/1-valued matrix that extracts off-diagonal elements, such as
\begin{equation}
    M_2 \text{vec} \left(\begin{bmatrix} 1 & 3 \\ 2 & 4\end{bmatrix}\right)= \begin{bmatrix} 2 \\ 3\end{bmatrix}.
\end{equation}

The polytope $P^C_2$ of continuous-time plants $(A, B)$ that can be positive-stabilized via \eqref{eq:stab_delta_c} under a state-feedback controller $K \in \scr$ with \iac{DLCLF} $\max(x./v)$ such that $Y = K X$ can be described by
\begin{subequations}
\label{eq:poly_stab_c}
\begin{align}
    G^C_2 &= \begin{bmatrix} v^T \otimes I_n & (Y\1_n)^T \otimes I_n \\
    -M_n( X \otimes I_n) & -M_n(Y^T \otimes I_n)\end{bmatrix} \\
    P^C_2 &= \left\{(A, B) \mid G^C_2 \begin{bmatrix} a \\ b \end{bmatrix} \leq \begin{bmatrix} -\eta \1_n \\ \0_{n(n-1)} \end{bmatrix}\right\}.
\end{align}
\end{subequations}

The top row of $G^C_2$ is the \ac{DLCLF} stabilization criterion, and the bottom row enforces that $AX+BY$ is Metzler.

The polytope $P^D_2$ of discrete-time plants $(A, B)$ positive-stabilized by $(K, Y)$ under the same conditions is
\begin{subequations}
\label{eq:poly_stab_d}
\begin{align}
    G^D_2 &= \begin{bmatrix} v^T \otimes I_n & (Y\1_n)^T \otimes I_n \\
     -X \otimes I_n & -Y^T \otimes I_n\end{bmatrix} \\
    P^D_2 &= \left\{(A, B) \mid G^D_2 \begin{bmatrix} a \\ b \end{bmatrix} \leq \begin{bmatrix} v -\eta \1_n \\ \0_{n^2} \end{bmatrix}\right\}.
\end{align}
\end{subequations}

\subsection{Stabilizing Programs using the Extended Farkas Lemma}
\label{sec:stab_program}

To unite notation, let $P_2$ be the appropriate stabilizing polytope for continuous-time $(P_2^C)$ or discrete-time $(P_2^D)$ from Section \ref{sec:stab_poly}. The number of constraints in the stabilizing polytope $P_2$ (length of $h_2$) is $q = n + n(n-1)$ for continuous-time and $q = n + n^2$ for discrete-time.
The polytope $P_2$ has a constraint matrix $G_2\in \R^{q \times n(n+m)}$ and vector $h_2 \in \R^{q}$ such that $P_2 = \{(A, B) \mid G_2 [a^T \ b^T]^T \leq h_2\}$. The entries in $G_2$ and $h_2$ are affinely-dependent on $(v, Y)$.

Problem \ref{prob:stab} may be expressed in the language of polytope-containment as,
\begin{prob}
\label{prob:poly_contain}
Find a vector $v \in \R^n_{> 0}$ and a matrix $Y \in \scr$ such that $P_1 \subseteq P_2$.
\end{prob}

\begin{thm}
\label{thm:poly_contain_lp}
Problem \ref{prob:poly_contain} (equivalent to \eqref{prob:stab}) has a solution iff the following \ac{LP} involving variables $(v, Y, Z)$ is feasible:
\begin{subequations}
\label{eq:lp_stab}
\begin{align}
    \find_{v, Y, Z} \qquad & Z G_1 = G_2(v, Y), \qquad Z h_1 \leq  h_2(v, Y) \label{eq:lp_stab_ext}\\
    & v-\eta \1_n \in \R_{\geq 0}^n, \ Y \in \scr, \ Z \in \R_{\geq 0} ^{q \times 2 n T}, \label{eq:lp_stab_var}
\end{align}
\end{subequations}
whereby the state-feedback gain $K \in \scr$ can be recovered by calculating $K = Y X^{-1}$.
\end{thm}
\begin{proof}
The \ac{LP} in \eqref{eq:lp_stab} is a direct application of the Extended Farkas Lemma \ref{lem:ext_farkas} to prove polytope containment $P_1 \subseteq P_2$.
\end{proof}

\subsection{Computational Complexity}

Table \ref{tab:num_cons} computes the number of inequality and equality constraints required to represent Program \eqref{eq:lp_stab_ext}. The number of equality constraints associated with $Y$ is set to 0 because zero-valued entries of $Y$ will be removed and will not be treated as scalar variables.
The \ac{LP} in \eqref{eq:lp_stab} has up to $n + mn + (2 n T) q$ scalar variables distributed into $(v, Y, Z)$, plus $q$ additional nonnegative slack variables required to represent the inequalities in constraint \eqref{eq:lp_stab_ext}.

\begin{table}[h]
    \centering
        \caption{Number of Inequality and Equality constraints in Program \eqref{eq:lp_stab}}
    \label{tab:num_cons}
    \begin{tabular}{c c c}
        & \# Ineq. & \# Eq.\\ \hline
        $v$ & $n$ & 0 \\
        $Y$ & $\leq mn$ & 0\\
        $Z$ & $(2 n T) q$ & 0\\
        \eqref{eq:lp_stab_ext} & $q$ & $q n(n+m)$
    \end{tabular}
\end{table}

In discrete-time with $q = n^2+n$ and no value-restrictions on $K$ ($Y \in \R^{m \times n}$), Program \eqref{eq:lp_stab_ext} will have $N = (2nT+1)(n^2 + n) + (2m + 1)n$ nonnegative scalar variables (representing $Y = Y^+ - Y^-$ where both $Y^+$ and $Y^-$ are nonnegative) and $(n^2+1)n(n+m)$ equality constraints.


The running-time of an Interior Point Method solver for \acp{LP} up to $\gamma$-optimality is approximately $O(N^{\omega+0.5} \abs{\log(1/\gamma)})$ \cite{wright1997primal}, where  $\omega$ is the matrix-multiplication constant. Our \ac{DDC} algorithm therefore has performance on the order of $(Tn^3)^{\omega+0.5} \sim n^{12.5}$. Significant gains in performance may be realized by noting that the matrices $(G_1, G_2)$ are sparse and are highly structured.

\begin{rmk}
\label{rmk:nonredundant_face}
The polytope $\Sigma_\dc$ may possess a large number of redundant faces. These half-space constraints may be removed to improve computational performance without affecting the description of $\Sigma_\dc$. Nonredundant faces may be discovered by linear programming over the polytope \cite{caron1989degenerate}. 
\end{rmk}

\begin{rmk}
An alternative approach is to perform vertex enumeration, in which relations \eqref{eq:stab_clean} hold at every vertex of $\Sigma_\dc$. The polytopes $\Sigma_\dc$ that are gathered as part of the data-acquisition process empirically have a number of vertices that scales exponentially with dimension, for which the face-based approach of the Extended Farkas Lemma is more favorable.
\end{rmk}

\begin{rmk}
This paper focused on the case of $L_\infty$-bounded noise. This set-containment framework will also be nonconservative when applied to other with other semidefinite-representable noise processes, such as when each column of the discrepency matrix $\mathbf{W}$ in \eqref{eq:discrepancy} has bounded $L_2$ norm. The Extended Farkas Lemma \ref{lem:ext_farkas} is a specific instance of a more general Robust Counterpart posed over a system of linear inequalities \cite[Theorem 1.3.14]{ben2009robust}.  In the $L_2$ case, each inequality constraint in the polytope in $P_2$ over the uncertain $(a, b)$ is replaced via a robust counterpart by $n(n+m)$ second-order-cone variables, $n(n+m)$ linear equality constraints, and one linear inequality constraint. This procedure is performed programmatically in \cite{lofberg2012robust} under the `duality' option.
\end{rmk}

\section{Peak-to-Peak Gain Regulation}
\label{sec:peaktopeak}
This section performs worst-case peak-to-peak (p2p) gain minimization using the Extended Farkas Lemma.

System \eqref{eq:sys} may be affected by an external noise process $\xi \in \R^e$ to form dynamics with a controlled output of $z \in \R^p$
\begin{subequations}
\label{eq:sys_xi}
\begin{align}
    \dx(t) &= A x (t) + B u(t) + E \xi(t) \\
    z(t) &= C x (t) + D u(t) + F \xi(t).
\end{align}
\end{subequations}
For a given set of parameters $(A, B, C, D, E, F)$ with $C \in \R^{p \times n}, \ D \in \R^{p \times m}, \ E \in \R^{n \times e}, \ F \in \R^{p \times e}$, this peak-to-peak gain may be computed by solving \iac{LP},
\begin{lem}[\cite{briat2013positive}]
\label{lem:p2p_clean}
There exists a state-feedback controller $u=Kx$ with $K, Y \in \scr$ and $v \in \R^n_{>0}$ such that peak-to-peak gain of \eqref{eq:sys_xi} is less than or equal to $\gamma \geq 0$ for continuous-time if
\begin{subequations}
\label{eq:p2p_c_clean}
\begin{align}
&-(A X + B Y)\1_n - E \1_e \in \R^{n}_{> 0} \label{eq:p2p_c_clean_dyn} \\
& \gamma \1_q - (C X + D Y)\1_n - F \1_e \in \R^{q}_{> 0} \\
& C X + D Y \in \R^{q \times n}_{\geq 0} \\
& AX + BY \ \textrm{is Metzler},
\end{align}
\end{subequations}
and for discrete-time if
\begin{subequations}
\label{eq:p2p_d_clean}
\begin{align}
&v -(A X + B Y)\1_n - E \1_e \in \R^{n}_{> 0} \label{eq:p2p_d_clean_dyn}\\
& \gamma \1_q - (C X + D Y)\1_n - F \1_e \in \R^{q}_{> 0} \\
& C X + D Y \in \R^{q \times n}_{\geq 0} \\
& AX + BY \in \R^{n \times n}_{\geq 0},
\end{align}
\end{subequations}
whereby the state-feedback gain can be recovered by $K = Y X^{-1}$.
\end{lem}



We aim to solve the following problem:
\begin{prob}
\label{prob:p2p}
Find a state-feedback controller $u=Kx$ with $K \in \scr$ to minimize the worst-case peak-to-peak gain $\xi \rightarrow z$ for any data-consistent plant $(A, B) \in  \Sigma_\dc$.
\end{prob}
\begin{rmk}
The $\epsilon$-corrupted data in $\dc$ is obtained when $\xi(t)=0$ at all time samples. It is further assumed that the matrices $(C, D, E, F)$ are all fixed and are known in advance.
\end{rmk}

Peak-to-peak polytopes for $(A, B)$ in \eqref{eq:p2p_c_clean} and \eqref{eq:p2p_d_clean} may be construct in a similar manner to the stabilizing polytopes $P_2$ in Section \ref{sec:stab_poly}. The right hand sides of these polytopes for \acp{CTS} and \acp{DTS} are,
\begin{align}
\label{eq:h_answer_p2p}
    h_2^{\textrm{p2p}:\, C} &=  \begin{bmatrix} -\eta \1_n - E \1_e\\ \0_{n(n-1)} \end{bmatrix} & h_2^{\textrm{p2p}:\, D} &=  \begin{bmatrix} v-\eta \1_n - E \1_e\\ \0_{n^2} \end{bmatrix}.
\end{align}

\begin{thm}
\label{thm:p2p}
Problem \ref{prob:p2p} has a solution iff the following \ac{LP} in variables $(v, Y, Z, \gamma)$ is feasible,
\begin{subequations}
\label{eq:p2p}
\begin{align}
    \gamma^* =&\min_{\gamma \in \R}  \qquad \gamma \label{eq:p2p_obj}\\
    & Z G_1 = G_2(v, Y), \qquad Z h_1 \leq h_2^{\textrm{p2p}}(v, Y) \label{eq:p2p_farkas}\\
   & (\gamma-\eta) \1_q - (C X + D Y)\1_n - F \1_e \in \R^{q}_{\geq 0} \\
& C X + D Y \in \R^{q \times n}_{\geq 0} \\
    & v-\eta \1_n \in \R_{\geq 0}^n, \ Y \in \scr, \ Z \in \R_{\geq 0} ^{q \times 2 n T}, \label{eq:lp_p2p_var}
\end{align}
whereby the p2p-minimizing state feedback gain $K \in \scr$ can be recovered by $K=Y X^{-1}$.
\end{subequations}
\end{thm}
\begin{proof}
The outer (peak-to-peak) polytope is $P_2^{\textrm{p2p}} = \{(A, B) \mid  G_2 [a^T \ b^T]^T \leq h_2^{\textrm{p2p}}\}$ for the appropriate continuous-time or discrete-time vector in  \eqref{eq:h_answer_p2p}, as constructed from conditions \eqref{eq:p2p_c_clean_dyn} or \eqref{eq:p2p_d_clean_dyn}. The Extended Farkas Lemma \ref{lem:ext_farkas} is then applied in \eqref{eq:p2p_farkas} to ensure that $\gamma$ is an upper bound for the peak-to-peak gain of all consistent systems. The objective in \eqref{eq:p2p_obj} reduces this gain as much as possible. The minimum is achieved because all constraints in \eqref{eq:p2p} are strict (due to the given tolerance $\eta > 0$).
\end{proof}


\section{Positive Switched Systems}
\label{sec:switched}

This section will extend the approach of Section \ref{sec:stab} to switched linear systems.

\subsection{Problem Setting}

The ground-truth switched system will be composed of $N_s$ subsystems with parameters $(A_s, B_s)$ for $s\in1..N_s$. A right-continuous switching function $S: [0, \infty) \rightarrow 1..N_s$ (\ac{CTS}) or $S: \N \rightarrow 1..N_s$ (\ac{DTS}) is used to define a switched system trajectory. Given a switching sequence $S$, the switched system trajectory  $x(t)$ satisfies 
\begin{align}
\label{eq:switched_traj}
    \delta x(t) &= A_{S(t)} x(t) + B_{S(t)} u_{S(t)}(t) & \forall t \geq 0.
\end{align}

System \eqref{eq:sys} is an instance of  \eqref{eq:switched_traj} with only a single subsystem $(A, B)$. We will assume in this  paper that $S$ may switch between subsystems arbitrarily often, and is not subject to dwell-time constraints. The matrices $B_{s} \in \R^{n \times m_s}$ are allowed to have different number of inputs $m_s$.

Data from \eqref{eq:data} is collected with the addition of a known switching sequence $\bs$ to form $\dc = (\bx, \bu, \xd, \bs)$: 
\begin{align}
\label{eq:data_switch}
    \begin{array}{ccllll}
        \bx & := & [x(0) & x(1) &  \ldots & x(T-1) ]   \\
        \bu & := & [u(0) & u(1) & \ldots & u(T-1)] \\ 
        \xd & := & [\delta x(0) & \delta x(1) &  \ldots & \delta x(T-1) ] \\
        \bs & := & [S(0) & S(1) &  \ldots & S(T-1) ].
    \end{array}
\end{align}

Our goal is to find \acp{DLCLF} and controllers to stabilize all possible systems \eqref{eq:switched_traj} consistent with \eqref{eq:data_switch}. We will discuss scenarios where the switching sequence $S(t)$ is unknown or is known.
 
\subsection{Data-Consistency Polytope}

We will collect together columns of $(\bx, \bu, \xd)$ that have the same switching subsystem. As an example, $(\bx^s, \bu^s, \xd^s)$ are the matrices formed from data in $\dc$ where $S(t) = s$. For notational convenience, we will use the abbreviations $a_s = \vvec{A_s}$ and $b_s = \vvec{B_s}$ for $s=1..N_s$. Given an $L_\infty$ noise bound $\epsilon$ for the discrepancy \eqref{eq:discrepancy}, we express the set of all $\dc$-consistent switched systems as the intersection of subsystem polytopes from \eqref{eq:data_single}
\begin{subequations}
\label{eq:data_switched}
\begin{align}
G^{\text{data}}_s &= \begin{bmatrix}(\bx^s)^T \otimes I_n & (\bu^s)^T \otimes I_n \end{bmatrix}\\
    \Sigma^{\text{data}}_{\dc, s} &=  \left\{ (A, B) \mid G^{\textrm{data}}_s \begin{bmatrix} a_s \\ b_s\end{bmatrix} \leq   \begin{bmatrix}\epsilon\1_{nT} +\vvec{\xd^s} \\ \epsilon\1_{nT} -\vvec{\xd^s}\end{bmatrix} \right\} \\
    \Sigma^{\text{data}}_{\dc} &= \bigcap_{s=1}^{N_s}\Sigma^{\text{data}}_{\dc, s}.
\end{align}
\end{subequations}

\subsection{Stabilizing Polytope}

This section will list stabilization tasks and their nominal \ac{DLCLF}-based criteria for stabilization. Each task description will conclude with a presentation of their \ac{DDC} stabilization polytope. 
The specific \ac{LP} of the form of Theorem \eqref{eq:lp_stab} will be skipped; it is enough to note that 
Extended Farkas Lemma \ref{lem:ext_farkas} will be used to compute stabilizing controllers based on the data-consistent polytope \eqref{eq:data_switched}  and the respective stabilizing polytope.

\subsubsection{Common DLCLF, Common Controller}
\label{sec:common_controller}

This task will require the number of inputs are the same among all subsystems $\exists m : m_s = m \forall s\in1..N_s$, and that the controller is unaware of the current switching sequence $S(t)$. 

\begin{thm}
\label{thm:positive_stable_common}
The switched system \eqref{eq:switched_traj} is positive-stabilized by the switching-independent state-feedback controller $u=Kx$ if there exists \iac{DLCLF} with $v > 0$, a matrix $X = \diag{v}$, a tolerance $\eta>0$, and a matrix $Y \in \R^{m \times n}$ such that $\forall s=1..N_s:$
\begin{subequations}
\begin{align}
    -(A_s X + B_s Y) \1_n -\eta \1_n& \in \R_{\geq 0}^n  & & \text{(\ac{CTS})} \label{eq:stab_delta_common_c} \\
    v-(A_s X + B_s Y) \1_n -\eta \1_n& \in \R_{\geq 0}^n & & \text{(\ac{DTS})} \label{eq:stab_delta_common_d}.
\end{align}
\end{subequations}
The control gain $K$ is recovered by $K=Y X^{-1}$.
\end{thm}
\begin{proof}
Switched-stabilization for a common \ac{DLCLF} is documented in \cite[Proposition 3.4]{blanchini2015switched} without applied control. The control term $Y$ is used as an extension of \cite{rami2007controller} (summarized in Theorem \ref{thm:stab_control}). The convergence of this \ac{DLCLF} scheme is established in \cite{pastravanu2014max}.
\end{proof}

Letting $P_2 = \{(A, B) \mid G_2[a^T b^T]^T \leq h_2\}$ be the single-system stabilization polytope (according to the notation in Section \ref{sec:stab_program}) with a common $(v, X, Y)$, the switched system stabilizing polytope is,
\begin{align}
    P_2^{\text{sw1}} &= \bigcup_{s=1}^{N_s} \{(A_s, B_s) \mid G_2[a_s^T, \ b_s^T]^T \leq h_2\}.
\end{align}

\subsubsection{Common DLCLF, Different Controller}
\label{sec:different_controller}
This task will allow the controller to access the switching sequence $S(t)$. Each subsystem $s$ is equipped with its own controller $K_s \in \R^{n \times m_s}$, in which the number of inputs $m_s$ is allowed to be different between subsystems. 

\begin{thm}
\label{thm:diff_sys}
System \eqref{eq:switched_traj} is positive-stabilized by the switching-dependent state-feedback controller $u(t)=K_{S(t)} x(t)$ if there exists \iac{DLCLF} with $v > 0$, a matrix $X = \diag{v}$, a tolerance $\eta>0$, and matrices $Y_s \in \R^{m_s \times n}$ such that $\forall s=1..N_s:$
\begin{subequations}
\begin{align}
    -(A_s X + B_s Y_s) \1_n -\eta \1_n& \in \R_{\geq 0}^n  & & \text{(\ac{CTS})} \label{eq:stab_delta_switch_c} \\
    v-(A_s X + B_s Y_s) \1_n -\eta \1_n& \in \R_{\geq 0}^n & & \text{(\ac{DTS})} \label{eq:stab_delta_switch_d}.
\end{align}
\end{subequations}
The subsystem control gains are recovered by $K_s = Y_s X^{-1}$ for $s=1..N_s$.
\end{thm}
\begin{proof}
The \ac{DLCLF} $\max(x./v)$ is common among all subsystems. Allowing the controller $K_s$ to change between subsystems allows more freedom to attempt solution for a positive-stabilizing controller.
\end{proof}

Let $G_2(v, X, Y)$ be the stabilizing polytope from the \ac{CTS} \eqref{eq:poly_stab_c} or the \ac{DTS} \eqref{eq:poly_stab_d} as appropriate. The switching-aware positive-stabilizing polytope $P^{\text{sw2}}_2$ is
\begin{align}
    P_2^{\text{sw2}} &= \bigcup_{s=1}^{N_s} \{(A_s, B_s) \mid G_2(v, X, Y_s)[a_s^T, \; b_s^T]^T \leq h_2\}.
\end{align}

\section{Positive Linear Parameter-Varying Systems}
\label{sec:lpv}


This section will perform data-driven positive-stabilzation for a class of \ac{LPV} systems. The presented approach is similar to subsystem-aware switched systems stabilization from Section \ref{sec:different_controller}.

\subsection{Problem Setting}

The \ac{LPV} framework involves parameters $\theta$ restricted to a known set $\Theta \subset \R^L$ that are measured on-line during operation. The general \ac{LPV} structure is
\begin{align}    \delta x &= A(\theta) x + B(\theta) u. \label{eq:lpv_general}
\intertext{We will focus on the \acf{LPVA} structure \cite{besselmannlofberg2012}, involving a set of matrices $\forall \ell: A_\ell \in \R^{n \times n}$ and a constant $B \in \R^{n \times m}$}
    \delta x &= \textstyle \left( \sum_{\ell=1}^L A_\ell \theta_\ell\right )x + B u. \label{eq:LPVA}
\end{align}
The open-loop system $\delta x = A(\theta) x$ is asymptotically stable if $\lim_{t\rightarrow \infty} x(t)=0$ for all possible parameter sequences $\theta(\cdot)$ taking values inside $\Theta$.

Data will be collected from \eqref{eq:lpv_general} with an $L_\infty$ noise bound of $\epsilon$. The observed data $\dc$ with $T$ records is 
\begin{align}
\label{eq:data_lpv}
    \begin{array}{ccllll}
        \bth & := & [\theta(0) & \theta(1) &  \ldots & \theta(T-1) ]\\
        \bx & := & [x(0) & x(1) &  \ldots & x(T-1) ]   \\
        \bu & := & [u(0) & u(1) & \ldots & u(T-1)] \\ 
        \xd & := & [\delta x(0) & \delta x(1) &  \ldots & \delta x(T-1) ]. \\
    \end{array}
\end{align}
The data in \eqref{eq:data_lpv} is collected into $\dc = (\bth, \bx, \bu, \xd)$.

Our goal is to find \iac{DLCLF} with parameter $v>0$ and a control policy $u(t) = K(\theta(t)) x(t)$ such that the closed-loop \ac{LPVA} system $\delta{x} = A(\theta) + BK(\theta)$ is stable and positive (Metzler for \ac{CTS} or Nonnegative for \ac{DTS}) for all $\theta \in \Theta$.

We will assume that there exists a finite $N_c \in \N$ and a bounded discrete set of points $\Omega = \{\omega_c\}_{c=1}^{N_c} \in \Theta$ such that $\Theta$ equals the convex hull of $\Omega$. We will refer to $\Omega$ as the set of `vertices' of $\Theta$. The controller will have knowledge of $\Omega$ and $\theta$ during operation.

\subsection{Data-Consistency Polytope}

Expression of the \ac{LPVA} data-consistency polytope will use the column-wise Khatri-Rao product for matrices $M_1 \in \R^{m \times n}, \ M_2 \in \R^{p \times n}$ (notated as $\kcol$)
\cite{khatri1968solutions}
\begin{equation}
    \label{eq:kron_col}
    M_1 \kcol M_2 = (\1_{p \times 1} \otimes M_1) \odot (M_2 \otimes \1_{m \times 1}).
\end{equation}

Data consistency of the \ac{LPVA} plant $\{A_\ell, B\}$ from \eqref{eq:LPVA} with $L_\infty$-norm error $\epsilon$ w.r.t. $\dc$ requires that
\begin{align}
    \forall t&\in0..T-1: & \norm{\dx(t) - \textstyle (\sum_{\ell=1}^L A_\ell \theta(t)_\ell) - B u(t)}_\infty \leq \epsilon. \label{eq:lpva_consistent_linf}
\end{align}

Define $a_\ell$ as $\vvec{A_\ell}$ for all $\ell=1..L$. Constraint \eqref{eq:lpva_consistent_linf} is equivalent to requiring that all $T$ columns of the following matrix $\bw$ have $L_\infty$-norm $\leq \epsilon$:
\begin{align}
\label{eq:lpv_noise}
    \bw &= \xd - \textstyle \left(\sum_{\ell=1}^L \bth_\ell \kcol A_\ell\right) \xn- B \bu.
\end{align} 

The data-consistent \ac{LPVA} polytope arising from \eqref{eq:lpv_noise} is
\begin{subequations}
\label{eq:lpv_consistent}
\begin{align}
    G^{\textrm{data}}_{LPV} &= \begin{bmatrix} (\bx \kcol \bth)^T \kron I_n & \bu^T \kron I_n\end{bmatrix} \\
    \Sigma^{\textrm{data}}_{LPV} &= \left\{ \{A_\ell, B\} \mid G^{\textrm{data}}_{LPV}  \begin{bmatrix} a_1^T, \ a_2^T, \ \ldots, \  a_L^T, b_T\end{bmatrix}^T \leq   \begin{bmatrix}\epsilon\1_{nT} +\vvec{\xd^s} \\ \epsilon\1_{nT} -\vvec{\xd^s}\end{bmatrix} \right\}
\end{align}
\end{subequations}

\subsection{Stabilizing Polytope}

Positive-Stabilization of consistent \ac{LPVA} will occur using a gain-scheduled controller \cite{rugh2000research} based on vertex-interpolation \cite{apkarian1995self}. 

The $\theta$-dependent control gain $K(\theta)$ will be constructed from a linear combination of controllers $\{K_c\}_{c=1}^{N_c}$. Each vertex $\omega_c \in \Omega$ has a corresponding vertex-controller $K_c \in \R^{m \times n}$ for every $c=1..N_c$. The control policy $u = K(\theta) x$ at a given $\theta$ will be found by first finding a feasible solution $\beta$ to the following \ac{LP}
\begin{subequations}
\label{eq:gain_sched}
\begin{align}
    \textrm{find} \ \beta & \in \R^{N_c}_+ &  \textstyle \sum_{c=1}^{N_c} \beta_c & = 1 &\textstyle   \sum_{c=1}^{N_c} \beta_c \omega_c & = \theta,
    \label{eq:interp_cert}
\end{align}
and subsequently returning the linear combination
\begin{align}
     K(\theta) &= \textstyle \sum_{c=1}^{N_c} \beta_c K_c & u &= K(\theta) x. \label{eq:gain_sched_out}
\end{align}
\end{subequations}

We define the vertex-plant $A_c$ corresponding to $\omega_c \in \Omega$ as
\begin{align}
    A_c &= \textstyle \sum_{\ell=1}^L \omega_{v \ell} A_\ell & & \forall v\in1..N_c
\end{align}

\begin{thm}
\label{thm:positive_stable}
The \ac{LPVA} system \eqref{eq:LPVA} is positive-stabilized by the gain-scheduled feedback controller $u=K(\theta) x$ from \eqref{eq:gain_sched_out} if there exists \iac{DLCLF} with $v > 0$, a matrix $X = \diag{v}$, a tolerance $\eta>0$, and matrices $Y_c \in \R^{m \times n}$ such that $\forall c=1..N_c:$
\begin{subequations}
\begin{align}
    -(A_c X + B Y_c) \1_n -\eta \1_n& \in \R_{\geq 0}^n  & & \text{(\ac{CTS})} \label{eq:stab_delta_lpva_c} \\
    v-(A_c X + B Y_c) \1_n -\eta \1_n& \in \R_{\geq 0}^n & & \text{(\ac{DTS})} \label{eq:stab_delta_lpva_d}.
\end{align}
The subsystem controllers $\{K_c\}_{c=1}^{N_c}$ may be recovered by $K_c = Y_c X^{-1}$.
\end{subequations}
\end{thm}
\begin{proof}
The \ac{LPVA} framework may be interpreted as switching in the sense of Theorem \ref{thm:diff_sys} between systems $(A_c, B)$ for $c\in1..N_c$. Compatibility between imposing positive-stabilization conditions on the vertices $\Omega$ and requiring the property to hold $\forall \theta \in \Theta$ is assured by Lemma 2.1 of \cite{miller2022lpvqmi}.
\end{proof}

The continuous-time \ac{LPVA} stabilization polytope from \eqref{eq:stab_delta_lpva_c} is
\begin{subequations}
\label{eq:poly_stab_lpva_c}
\begin{align}
    G^C_{2c}  &= \begin{bmatrix} \omega_c^T \kron (v^T \otimes I_n) & (Y\1_n)^T \otimes I_n \\
    -M_n( \omega_c^T \kron (X \otimes I_n) )& -M_n(Y^T \otimes I_n)\end{bmatrix} \\
    P^C_{2 LPV} &= \left\{(\{A_\ell\}, B) \mid \ \forall c=1..N_c: \ G^C_{2c} \begin{bmatrix} a_1^T, a_2^T, \ldots, a_L^T,  b^T \end{bmatrix}^T \leq \begin{bmatrix} -\eta \1_n \\ \0_{n(n-1)} \end{bmatrix} \right\}.
\end{align}
\end{subequations}

The discrete-time \ac{LPVA} stabilization polytope from \eqref{eq:stab_delta_lpva_d} is
\begin{subequations}
\label{eq:poly_stab_lpva_d}
\begin{align}
    G^D_{2c}  &= \begin{bmatrix} \omega_c^T \kron (v^T \otimes I_n) & (Y\1_n)^T \otimes I_n \\
    -( \omega_c^T \kron (X \otimes I_n) )& -(Y^T \otimes I_n)\end{bmatrix} \\
    P^D_{2 LPV} &= \left\{(\{A_\ell\}, B) \mid \forall c=1..N_c: \ G^D_{2c} \begin{bmatrix} a_1^T, a_2^T, \ldots, a_L^T,  b^T \end{bmatrix}^T \leq \begin{bmatrix} v-\eta \1_n \\ \0_{n^2} \end{bmatrix} \ \right\}.
\end{align}
\end{subequations}

The stabilizing polytopes $P^C_{2 LPV}$ and $P^D_{2 LPV}$ can be used in conjunction with the $\dc$-consistency polytope $\Sigma^{\textrm{data}}_{LPV}$ to form \ac{LPVA} \ac{DDC} \acp{LP} by the Extended Farkas Lemma.
\section{Numerical Examples}

\label{sec:examples}

MATLAB 2021a Code to reproduce the experiments is available at \url{https://github.com/jarmill/data_driven_pos}, and includes Mosek \cite{mosek92} and YALMIP \cite{lofberg2004yalmip} dependencies. All provided experiments have parameters of $\eta = 10^{-3}$ and $\epsilon = 0.1$.
 
\subsection{Continuous-Time Stabilization}

The ground-truth continuous-time system in this example has $n=3$ inputs and $m=2$ outputs
\begin{align}
\label{eq:cont_stab_sys_true}
    A &= \begin{bmatrix} -0.55 & 0.3 & 0.65\\ 0.06 & -1.35 & 0.25 \\ 0.1 & 0.15 & 0.4   \end{bmatrix}& B = \begin{bmatrix} 0.18 & 0.08 \\ 0.47 & 0.25 \\0.07 & 0.95 \end{bmatrix}.
\end{align}

System \eqref{eq:cont_stab_sys_true} is internally positive but is open-loop unstable (poles of $0.4907, -0.6055, -1.3851$). The stabilization task in \eqref{eq:lp_stab} with $T=5$ and an additional normalization constraint that $\1_n^T v = 1$ results in
\begin{subequations}
 \begin{align}
     v &= \begin{bmatrix}0.5570 & 0.1401& 0.3029\end{bmatrix}^T \label{eq:example_cont_stab_v}\\
    K &= \begin{bmatrix} 0.0279 &    -0.2660 &  0.5041 \\
    0.0107  &  -0.0222 & -0.8650\end{bmatrix}. \label{eq:example_cont_stab_K}
 \end{align}
\end{subequations}

\begin{figure}[h]
    \centering
    \includegraphics[width=0.5\linewidth]{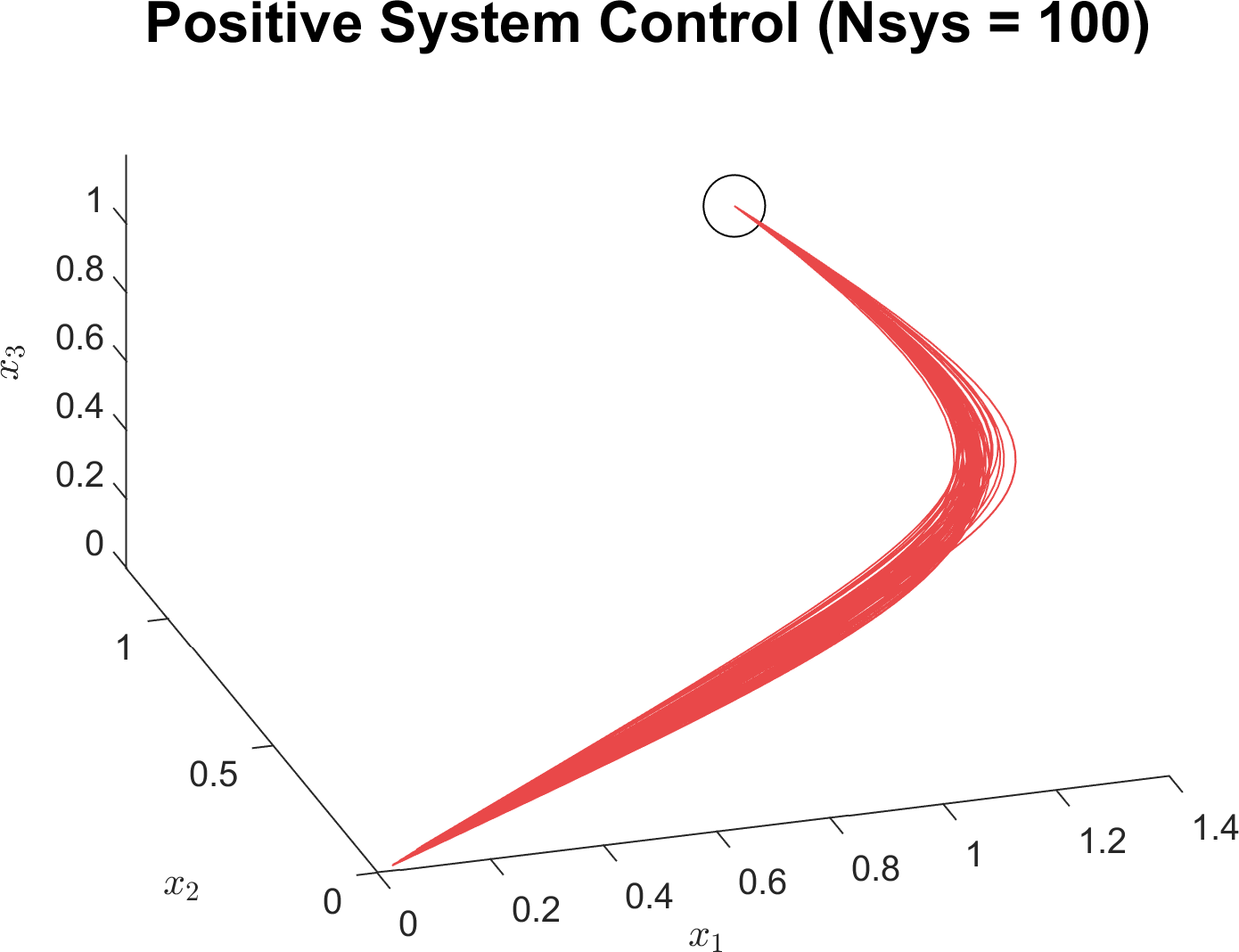}
    \caption{Application of the controller $u=Kx$ from \eqref{eq:example_cont_stab_K} to positively-stabilize $100$ consistent systems in $\Sigma_D$.}
    \label{fig:cont_stab}
\end{figure}

Figure \ref{fig:cont_stab} visualizes 100 controlled trajectories (red curves) starting from $x(0)=[1;1;1]$ (black circle). Each trajectory follows $\dot{x}(t) = (A+ BK)x(t)$ in the times $t \in [0, 20],$ where the plants $(A, B)$ are randomly sampled from $\Sigma_D$ and $K$ is the controller in \eqref{eq:example_cont_stab_K}.

\begin{figure}[h]
    \centering
    \includegraphics[width=0.75\linewidth]{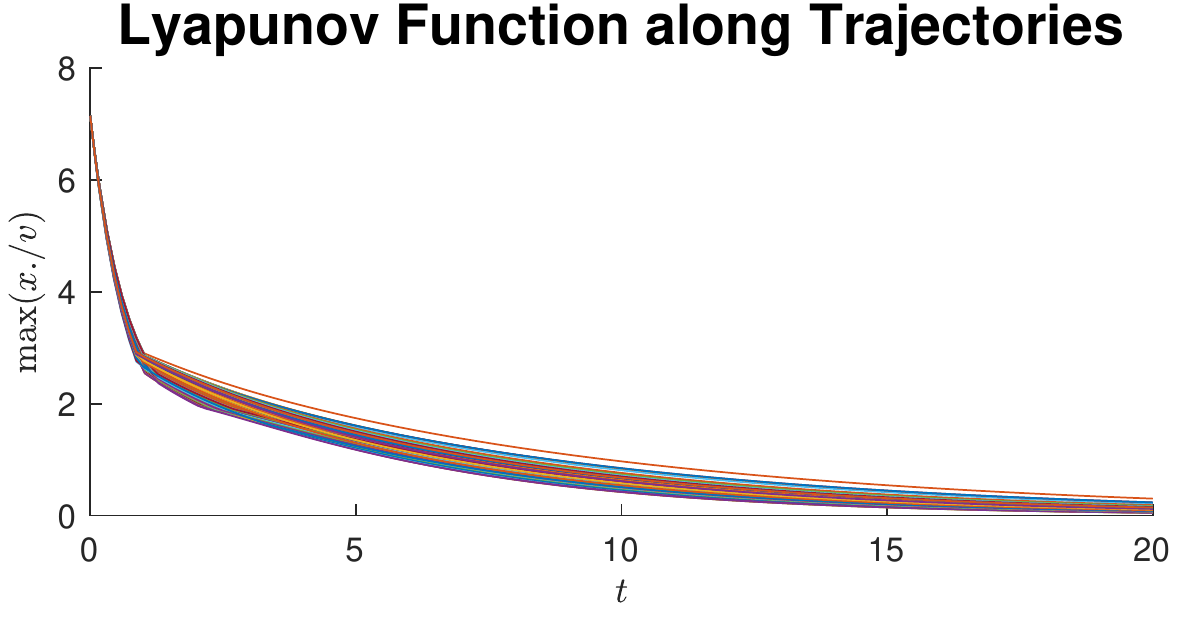}
    \caption{\ac{DLCLF} along the 100 trajectories.}
    \label{fig:cont_max_lyap}
\end{figure}

Figure \ref{fig:cont_max_lyap} plots values of the Lyapunov function $\max(x./v)$ (for the $v$ in \eqref{eq:example_cont_stab_v}) along the 100 systems in \ref{fig:cont_stab}.

\subsection{Discrete-Time Stabilization}

This example involves a discrete-time system with $n=5$ states and $m=3$ inputs. The ground-truth system is internally positive, and is unstable with poles of $1.3094, \ -0.1218 \pm 0.0992 \mathbf{j}, \-0.1201 \pm 0.1108 \mathbf{j}$.  With $T=60$ observations the following \ac{DLCLF} and stabilizing controller is recovered
\begin{align}
v &= \begin{bmatrix}0.2076& 0.1212&   0.2651&   0.2516 & 0\end{bmatrix}^T     \\
K &= \begin{bmatrix}0.0483 &   0.0088  & -0.1326  & -0.0188  & -0.4273\\
   -0.3243  &  0.0115  &  0.0299 &  -0.2980 &   0.0337\\
    0.1601 &   0.0749  & -0.5962 &  -0.3537 &  -0.2194 \end{bmatrix}. \nonumber
\end{align}

 
It is now desired to obtain a stabilizing controller for all consistent plants that obeys the sign pattern
\begin{subequations}
\begin{align}
    \scr &=  \begin{bmatrix}\odot &   \odot  & \odot  & \odot  & \ominus\\
   \odot  &  \odot  &  \circledast & \odot &   \oplus\\
    \odot &   \odot  & \odot &  \circledast &  \circledast \end{bmatrix}
\intertext{Such  \iac{DLCLF} certificate  and controller is}
v &= \begin{bmatrix} 0.2147&   0.1259&    0.2448&     0.2516&     0.1630\end{bmatrix}^T\\
K &= \begin{bmatrix}0&   0  & 0  & 0  & -0.6853\\
   0  &  0  &  -0.3206 & 0 &   0.1206\\
    0 &   0  & 0 &  -0.5604 &  -0.3317 \end{bmatrix}.
\end{align}

\end{subequations}
\subsection{Continuous-Time Peak-to-Peak}

The following ground-truth positive-stable continuous-time system has $n=3$ inputs and $m=2$ outputs
\begin{align}
\label{eq:p2p_sys_true}
    A &= \begin{bmatrix} -0.2 & 0.2 & 0.2\\ 0.4 & -0.7 & 0.2 \\ 0 & 0.8 & -3   \end{bmatrix}& B = \begin{bmatrix} -0.4 & 0.5 \\ 0.2 & -0.8 \\ -1 & 2 \end{bmatrix}.
\end{align}

This system has  $e=2$ external input channels and $p=5$ controlled outputs with
\begin{align}
\label{eq:p2p_param}
    C &= \begin{bmatrix} I_3 \\ \0_{2 \times 3}
    \end{bmatrix}, & D &= \begin{bmatrix} \0_{3 \times 2} \\ I_2
    \end{bmatrix}, & E  &= \begin{bmatrix} I_2 \\ \0_{1 \times 2}
    \end{bmatrix}, & F  &= \0_{5 \times 2}.
\end{align}
 
The peak-to-peak gain of the ground-truth \eqref{eq:p2p_sys_true} under the parameters in \eqref{eq:p2p_param} when uncontrolled $(K=\0_{2 \times 3})$ is $\gamma^* = 32.178$. Lemma \ref{lem:p2p_clean} synthesizes a controller for the ground-truth system resulting in a gain of $\gamma^* = 3.742$. The constraint $CX + DY \in \R^{q \times n}_{\geq0}$ with the values in \eqref{eq:p2p_param} imposes that all elements of $Y$ and $K$ are nonnegative $(\oplus)$.

Table \ref{tab:p2p_result} collects the worst-case peak-to-peak gains obtained by \eqref{eq:p2p} as a function of the number of samples $T$. These gains decrease as $T$ increases and the consistency set $\Sigma_\dc$ shrinks. The top row of \eqref{tab:p2p_result} incorporates the prior knowledge that the ground-truth $A$ from \eqref{eq:p2p_sys_true} is Metzler when constructing the polytope $\Sigma_D$. The bottom row does not impose this positivity (Metzler) prior on $A$, and therefore yields peak to peak bounds that are always greater than or equal to the Metzler-imposed bounds.

 \begin{table}[h]
     \centering 
     \caption{Worst-case peak-to-peak gain $\gamma^*$ computed by \eqref{eq:p2p} decreases as the number of samples $T$ increases}
     \begin{tabular}{rccccc}
     $T$& 20 & 30 & 50 & 80 & 120 \\ \hline
          $A$ Metzler & 6.4539 & 5.0182 & 4.4967 & 4.0619 & 4.0028 \\
          No Prior& 6.4823 & 5.0719 & 4.5292 & 4.0659 & 4.0029
     \end{tabular}
     \label{tab:p2p_result}
 \end{table}
 
The system with $T=50$ and a Metzler-prior on $A$ has a worst-case peak-to-peak gain of $\gamma^*=4.4967$ and solution outputs of
\begin{subequations}
 \begin{align}
     v &= \begin{bmatrix}4.4967 &4.2021 & 0.4303\end{bmatrix}^T \\
    K &= \begin{bmatrix} 0.5095 &    0.4765  &  0.4727 \\
    0.2587  &  0 & 0\end{bmatrix}.
 \end{align}
\end{subequations}
 
The polytope $\Sigma_D$ under the Metzler-prior has $2nT + (n^2-n) = 300 + 6 = 306$ faces and 308,672 vertices, of which 62 faces are nonredundant (see Remark \ref{rmk:nonredundant_face}). The nonnegative Farkas matrix is $Z \in \R^{9 \times 62}_{\geq 0}$.


\subsection{Switched System Control}

This example will involve a continuous-time system with $n=3$ inputs, $m=2$ outputs, and $N = 2$ subsystems
\begin{align}
    A_1^{\textrm{true}} &= \begin{bmatrix} -0.55 & 0.3 & 0.65\\ 0.06 & -1.35 & 0.25 \\ 0.1 & 0.15 & 0.4   \end{bmatrix}& B_1^{\textrm{true}} &= \begin{bmatrix} 0.18 & 0.08 \\ 0.47 & 0.25 \\0.07 & 0.95 \end{bmatrix} \nonumber\\
    A_2^{\textrm{true}} &= \begin{bmatrix} 0.1& 0.1 & 0.1\\ 0.1 & -1.9 & 0.15 \\ 0.1 & 0.1 & 0.6   \end{bmatrix}& B_2^{\textrm{true}} &= \begin{bmatrix} 1 & 0\\ 0 & 0 \\0 & 1 \end{bmatrix} \label{eq:three_state_common}. 
\end{align}

A set of $T=55$ observations of system \eqref{eq:three_state_common}, 28 of which in 
are $s=1$ and the remaining 27 in $s=2$. The polytope $P_1^{c}$ has $N n(n+m) = 30$ dimensions, 108 nonredundant faces, and 246 redundant faces. 

\begin{figure}[ht]
    \centering
      \includegraphics[width=0.85\linewidth]{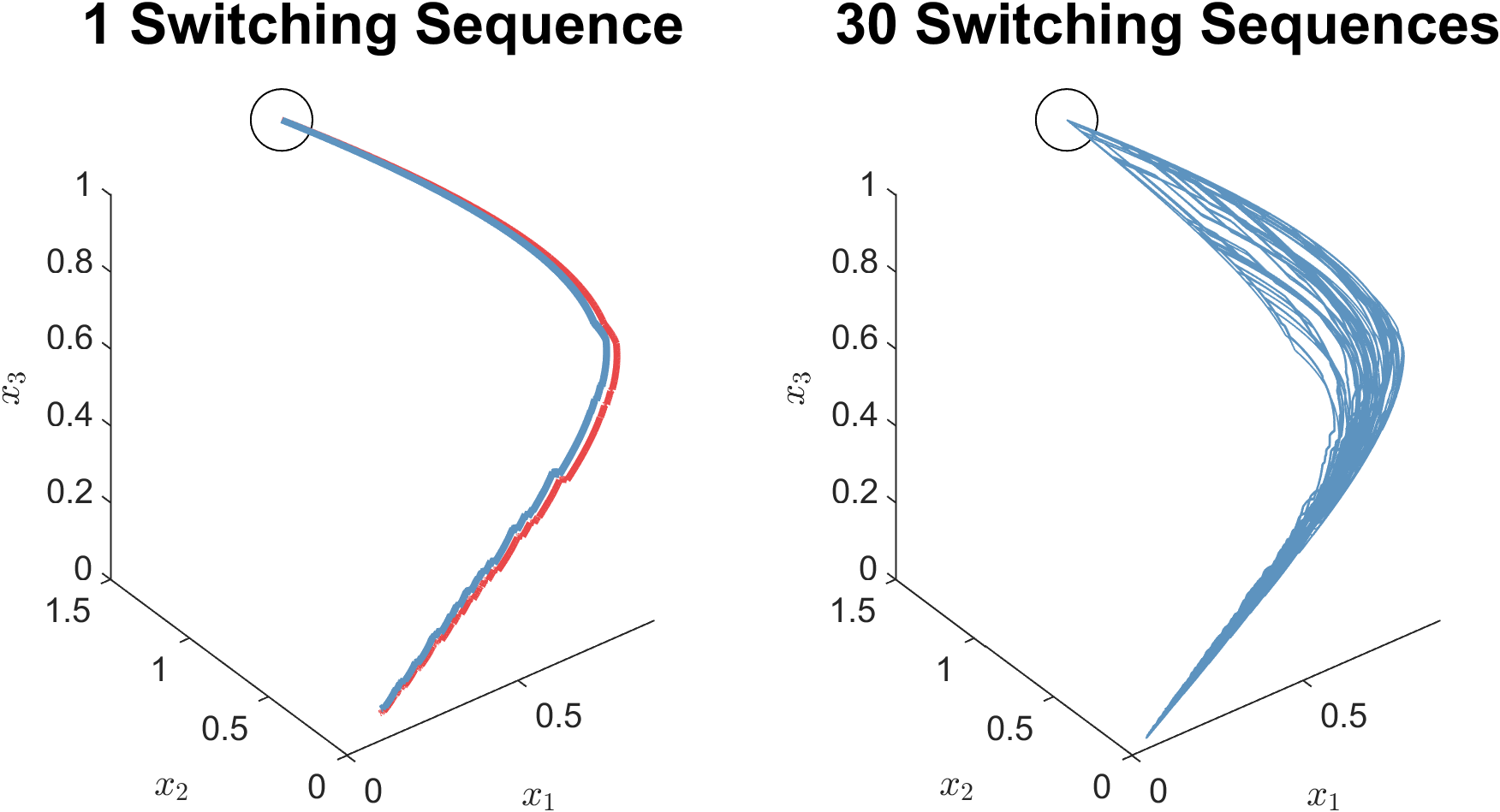}
    \caption{Controlled switched trajectories using the gain in \eqref{eq:three_state_common_K} }
    \label{fig:cont_common}
\end{figure}

The recovered controller (and \ac{DLCLF} vector) that simultaneously stabilizes both systems in \eqref{eq:three_state_common} (Section \ref{sec:common_controller}) are
\begin{subequations}
\begin{align}
    v &= \begin{bmatrix}
        0.4989 & 0.0572 & 0.4439
    \end{bmatrix} \\
    K &= \begin{bmatrix}        
   -0.1390 &  -0.0860  & -0.0663\\
    0.0362 &  -0.0810  & -0.8146
    \end{bmatrix} \label{eq:three_state_common_K}.
\end{align}
\end{subequations}
Figure \ref{fig:cont_common} plots controlled trajectories of the system in \eqref{eq:three_state_common} with the gain in \eqref{eq:three_state_common_K} starting from $x(0) = [0.5, 1.5, 1]$. The switching time of each trajectory to a new subsystem is exponentially distributed with a mean of 0.3 time units. The red trajectory on the left subplot highlights the ground truth system in \eqref{eq:three_state_common}, and the other blue curves are trajectories of 15 subsystems inside $P_1^\textrm{c}$ when the identical switching sequence is applied. The right subplot overlays trajectories of 30  switching sequences. 

Figure \ref{fig:cont_common} is generated with $T=55$ datapoints. When only $T=20$ observations are collected, it is infeasible to find a common \ac{DLCLF} and controller. However, a common \ac{DLCLF} and a pair of subsystem controllers that can stabilize both systems in \eqref{eq:three_state_common} (Section \ref{sec:different_controller}) are
 
\begin{subequations}
\begin{align}
    v &= \begin{bmatrix}
        0.5423 & 0.1327 & 0.3250
    \end{bmatrix} \\
    K_1 &= \begin{bmatrix}        
    0.0444 &  -0.3097  &  0.4207 \\
   -0.0010 &   0.2910  & -1.0869
    \end{bmatrix} \\
    K_2 &= \begin{bmatrix}        
   -0.4223  &  0.1510  &  0.1520 \\
    0.0059  & -0.0171  & -0.9607
    \end{bmatrix}\label{eq:three_state_diff_K}.
\end{align}
\end{subequations}

\begin{figure}[ht]
    \centering
      \includegraphics[width=0.85\linewidth]{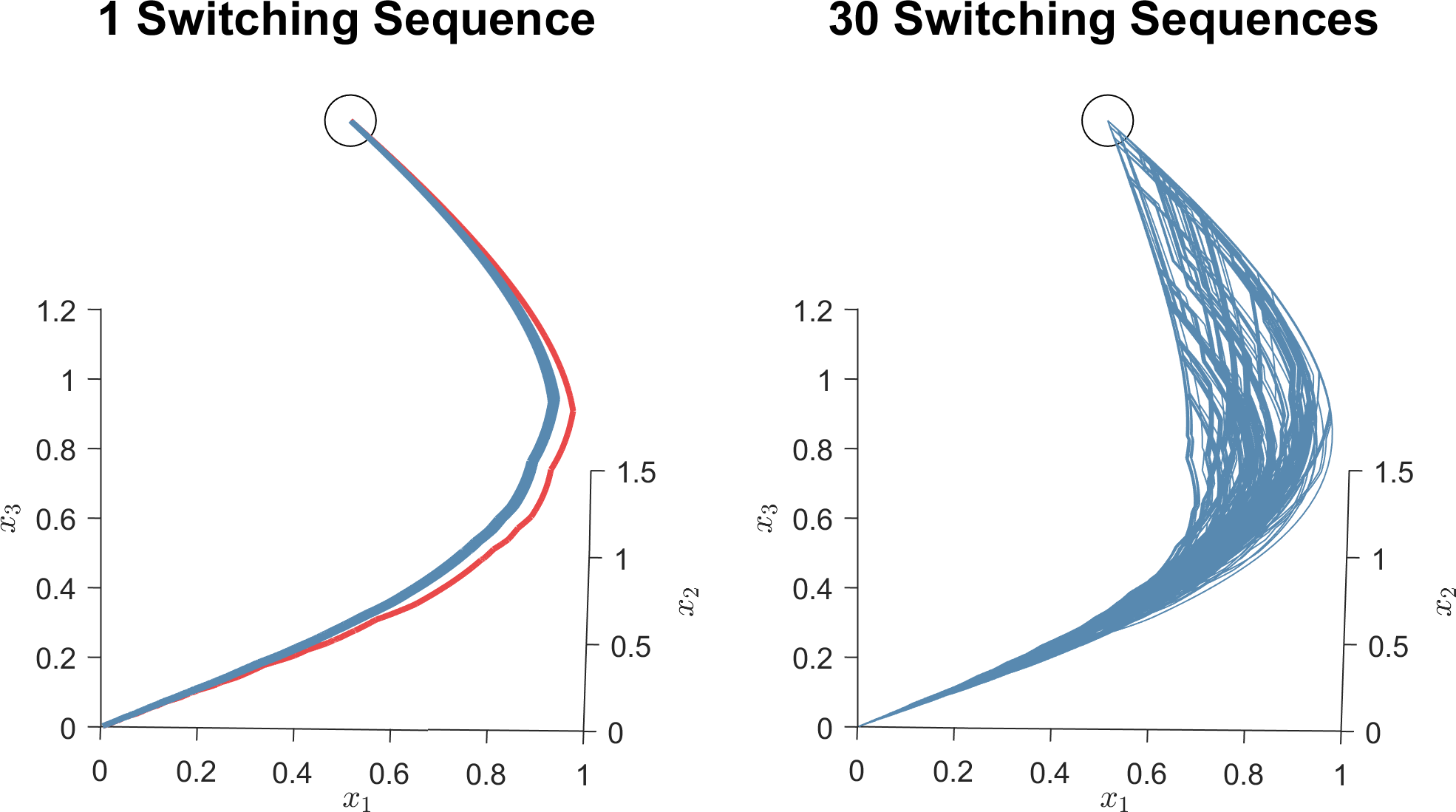}
    \caption{Controlled switched trajectories using the subsystem-dependent gains in \eqref{eq:three_state_diff_K} }
    \label{fig:cont_diff}
\end{figure}

Figure \ref{fig:cont_diff} plots switching-aware controlled trajectories of \eqref{eq:three_state_common}   based on \eqref{eq:three_state_diff_K}.


\subsection{LPV System}

The considered ground truth continuous-time system with a parameter set
$\Theta = \{1\} \times  [-1, 1] \times [-0.5,0.9]$ is
\begin{align}
    A_1^{\textrm{true}} &= \begin{bmatrix} -0.9190  & 0.5555 \\ 0.4936 & -0.5761 \end{bmatrix}, \nonumber &     A_2^{\textrm{true}} &= \begin{bmatrix} -1.2653  & 0.0574 \\ 0.2981 & 0.2455 \end{bmatrix} \nonumber\\
    A_3^{\textrm{true}} &= \begin{bmatrix} 0.9328  & 0.5702 \\ 0.0636 & -1.0487 \end{bmatrix}, & B^{\textrm{true}} &= \begin{bmatrix} 0.4570 & 0.2828 \\ 0.2115 & 0.8863 \end{bmatrix}.\label{eq:two_state_lpv}
\end{align}

This internally positive system has parameters $n = 2, m = 2$ and $L = 3$.
The plant matrix corresponding to the vertex $\omega = [1, -1, 0.9] \in \Theta$ is unstable, because $A_1^{\textrm{true}} - A_2^{\textrm{true}} + A_3^{\textrm{true}}$ has a positive eigenvalue of $1.2700$.

Data with a horizon of $T=10$ (10 state-input-transition tuples) was collected with $\epsilon =0.1$. The polytope $\Sigma_{LPV}^{\textrm{data}}$ from \eqref{eq:lpv_consistent} has $n(Ln + m) = 16$ dimensions,  42 nonredundant faces, 8 redundant faces, and 607590 vertices.

The following continuous-time vertex controllers positive-stabilize all consistent $\dc$-plants by \eqref{eq:stab_delta_lpva_c} using the stabilizing polytope \eqref{eq:poly_stab_lpva_c}.

\begin{align}
    K_{(1, -1, -0.5)} &= \begin{bmatrix}
  -14.2950  &  9.9057 \\
    6.2326  & -5.8745
    \end{bmatrix}   \nonumber\\  
    K_{(1, -1, 0.9)} &= \begin{bmatrix}
     -20.8043  &  9.7975 \\
    8.5165  & -6.9282
    \end{bmatrix} \nonumber \\
        K_{(1, 1, -0.5)} &= \begin{bmatrix}
   -6.9969 &   6.5542 \\
    2.5340 &  -4.3078
    \end{bmatrix} \label{eq:two_state_k} \\      K_{(1, 1, 0.9)} &= \begin{bmatrix}
   -5.2847  &  2.3813 \\
    2.0480  & -2.3017
    \end{bmatrix}. \nonumber
\end{align}

These controllers were synthesized under the prior knowledge $(\Sigma_{\textrm{prior}})$ that each ground-truth $\{A_\ell\}$ is Metzler and $B$ is nonnegative.

The controllers in \eqref{eq:two_state_k} have an associated common \ac{DLCLF}
\begin{align}
   V(x) &= \max(x_1/0.4482, x_2/0.5518). 
\end{align}

Figure \ref{fig:param_seq} plots system trajectories starting from the black-circle point $x(0) = [0.5, 1.5]$. Parameter values in $\theta$ are drawn uniformly from the box $\Theta$, and change values at mean-$0.05$ exponentially distributed switching times. The top plot of Figure \ref{fig:param_seq} plots controlled trajectory execution for the ground truth in \eqref{eq:two_state_lpv} (red) as well as 15 other systems randomly drawn from $P_1^{\textrm{c,LPV}}$. The bottom plot displays controlled trajectories arising from 30 parameter-switching sequences for each of the 16 sampled systems.

\begin{figure}[ht]
    \centering
      \includegraphics[width=0.6\linewidth]{img/sys_lpv_cont.png}
    \caption{Controlled \ac{LPV} trajectories using the gains in \eqref{eq:two_state_k} }
    \label{fig:param_seq}
\end{figure}

\section{Conclusion}

\label{sec:conclusion}

This paper presented an LP-based algorithm (Theorem \ref{thm:poly_contain_lp}) to perform data-driven stabilizing control of positive linear systems. The state-feedback controller $K$ stabilizes all possible systems in the $L_\infty$-norm bounded consistency set $\Sigma_\dc$, as certified by a common \ac{DLCLF} function $V(x) = \max(x./v)$ and the Extended Farkas Lemma. There is no conservativeness in such a design: Equation \eqref{eq:lp_stab} will find a controller iff there exists such a linear copositive Lyapunov function across all consistent systems. This framework can also be used to perform data-driven worst-case peak-to-peak gain minimization using Equation \eqref{eq:p2p}. The \ac{LP} framework is extended to positive switched and positive \ac{LPV} systems.

Future work includes applying these techniques to monotone systems and systems with dwell time. Other aspects include applying Lyapunov-Metzler inequalities to perform stabilization when the controller is able to select the switching sequence $S$ \cite{blanchini2015switched, bianchi2022data}. 




\bibliographystyle{IEEEtran}
\bibliography{references.bib}

\end{document}